\documentclass[10pt]{article}
\usepackage[letterpaper, margin=2.5cm]{geometry}

\usepackage{color, amsmath,amssymb, amsfonts, amstext,amsthm, latexsym}
\usepackage{epsfig}
\usepackage{longtable}
\usepackage{graphicx}
\usepackage{subfigure}
\usepackage{hyperref}
\usepackage{indentfirst}
\usepackage{url}

\renewcommand{\phi}{\varphi}

\begin{document}
\title{On the abrupt change of the maximum likelihood state in a simplified stochastic thermohaline circulation system}
 \author{\small{\it{Fang Yang$^{1,2}$, Xu Sun$^{1,2}$\footnote{Corresponding author: xsun15@iit.edu}, Jinqiao Duan$^{3}$}}\\
	\\ \small{$^1$ School of Mathematics and Statistics, Huazhong University of Science and Technology,
		Wuhan 430074, Hubei, China} \\
	\small{$^2$ Center for Mathematical Sciences, Huazhong University of Science and Technology,
		Wuhan 430074, Hubei, China}\\
	\small{$^3$Department of Applied Mathematics, Illinois Institute of Technology, Chicago, IL 60616, USA.}\\
}
\bigskip

\date{\today}
\maketitle

\pagestyle{plain}

Abstract: The maximum likelihood state for a simplified stochastic thermohaline circulation model is investigated. It is shown that a jump occurs for the maximum likelihood state during transitions between two metastable states. The jump helps to explain the mechanism of the abrupt change in the climate systems as revealed by various climate records.


\bigskip 

\textbf{It is well known that strong interaction exists between the thermohaline circulation system and the climate system. To understand the mechanism of the abrupt climate change, various methods have been proposed to estimate the likelihood of abrupt change in the thermohaline circulation system. We study the maximum likelihood state of a simplified stochastic thermohaline circulation system. It is shown that jumps  occur along transition pathways of the maximum likelihood state. The jump indicates a dramatic change of the salinity in the thermohaline circulation system. This observation is helpful to understand the mechanism of the abrupt change in the climate system.}

\section{Introduction}
 Driven by surface buoyancy forcing and pressure differences in the deep ocean, the thermohaline circulation depends heavily  on the temperature and salinity of seawater. The thermohaline circulation plays an important role in maintaining the overall energy balance of the Earth \cite{Weaver1992,Monahan2002}, and even a small change in the circulation may have a huge impact on the global climate.

 Various climate records reveal that the abrupt climate change is highly related to the transition between the equilibrium states of the thermohaline circulation \cite{Velez2001}.   Note that in deterministic systems,  transitions between equilibrium states  would not happen without external excitations.  Transitions from one equilibrium state to another are often explained from the viewpoint of stochastic models \cite{Cessi1994,Timmermann2000,Lorenzo2008}. It is shown \cite{Velez2001}  that the thermohaline circulation is  a nonlinear bistable system subjected to fluctuating environment, and the stochastic freshwater flux can push the thermohaline circulation system to transit from one stable state to another. Some approaches based on past experiences, such as observed climate data \cite{Knight2005,Zhang2010} or model simulation \cite{Averyanova2017,Rahmstorf2002}, have been proposed to estimate the likelihood of abrupt change in the thermohaline circulation system.
 
Hoping to provide some insights into the mechanism of abrupt change of the state in the thermohaline circulation system,   we  estimate the maximum likelihood state of the system and investigate how it  evolves during the transition process from one equilibrium point to the other.

This paper is organized as follows. In section 2,  the stochastic thermohaline circulation  model is introduced. The maximum likelihood state transition is investigated with the help of simulation results in section 3. Section 4 is the conclusion.

\bigskip
\section{A stochastic thermohaline circulation  model}

One hemisphere thermohaline circulation can be modeled by two vessels of saltwater solution that are connected by advective flow from below and exchange thermal energy and salinity with each other and with the surrounding  atmosphere. One vessel represents the high latitudes or north pole with temperature $T_p$ and salinity $S_p$, the other vessel stands for the low latitudes or equatorial region having temperature $T_e$ and salinity $S_e$. Each of the vessels is well mixed and has uniform salinity and temperature. But the salinity and temperature between the vessels may not be the same. 

The temperature and salinity differences between equatorial and pole regions, defined as  $\Delta T \equiv T_e-T_p$ and $\Delta S \equiv S_e-S_p$, respectively, are shown to evolve with time by \cite{Cessi1994}
\begin{equation} \label{e7}
\begin{split} 
	\dfrac{{\rm d} \Delta T}{dt } &=-\dfrac{1}{t_r}\left(\Delta T -\theta\right)-Q(\Delta \rho) \Delta T,\\
	\dfrac{{\rm d}  \Delta S}{dt } &= \dfrac{F}{H}S_0-Q(\Delta \rho)\Delta S,
	\end{split}
\end{equation}
where $\Delta \rho$ is the density difference expressed as
\begin{align}\label{e5}
\Delta \rho=\alpha_{S}(S_p-S_e)-\alpha_{T}(T_p-T_e),
\end{align}
and $Q(\Delta \rho)$ is  the exchange function expressed as
\begin{align}\label{e6}
Q(\Delta \rho)= \dfrac{1}{t_d}+\dfrac{q}{V}\left( {\Delta \rho}\right)^{2}.
\end{align}

Following Cessi \cite{Cessi1994},  the values of the dimensional parameters for the  thermohaline circulation model \eqref{e7} are shown in Table 1.  
With the scaling on temperature difference, salinity difference and time, $x\equiv \dfrac{\Delta T}{\theta}, y\equiv \dfrac{\alpha_S \Delta S}{\alpha_{T} \theta}, t'\equiv \dfrac{t}{t_d}$,  \eqref{e7} can be rewritten as the following non-dimensional system,
\begin{align} \label{e8}
\begin{split}
	\dot{x}&=-\alpha(x-1)-x\left(1+\mu^{2}(x-y)^2 \right),\\
	\dot{y}&=\bar {F}-y\left(1+\mu^2 (x-y)^2\right),
\end{split}
\end{align}
where
\begin{align}\label{e10}
\bar{F}=\dfrac{\alpha_SS_0t_d}{\alpha_{T}\theta H}F,~~~\alpha=\dfrac{t_d}{t_r}, ~~~ \mu^2 =\dfrac{t_d}{t_a}, ~~~t_a=\dfrac{q(\alpha_{T}\theta)^2}{V},
\end{align}
 the non-dimensional parameter $\alpha$ measures temperature restoring tensity, $\mu^2$  the ratio of the diffusion time scale $t_d$ and the temperature relaxation time scale $t_a$, and $\bar{F}$  the non-dimensional freshwater flux. 

\begin{table}[h]
	\caption{\label{tab:table1}Parameters of the thermohaline circulation  model given by \eqref{e7}}
	\centering
	\begin{tabular}{llll}
		\hline
		Parameter&Description& Value&Dimension\\
		\hline
		$T_0$& reference  temperature& $5$& $^{\rm o}{\rm C}$\\
		$S_0$& reference salinity & $35$& ${\rm psu^{-1}}$\\
		$\rho_0$& reference density & $1029$& ${\rm Kg~m^{-3}}$\\
		$t_d$                & the diffusion time& 219& ${\rm yeays}$\\
		$t_r$ & temperature relaxation time & 25&${\rm days}$\\
	    $t_a$ & advection time scale& 35&${\rm years}$\\
		$q$ &transport coefficient &$8.84\times 10^{11}$&${\rm m^3s^{-1}}$\\
		$\alpha_{S}$ & salinity coefficient& $7.5$ $ \times 10^{-4}$ &${\rm psu^{-1}}$\\
		$\alpha_{T}$ & thermal expansion coefficient  & $1.7\times 10^{-4}$&$^{\rm o}{\rm C^{-1}}$\\
		$\theta$ & meridional temperature difference &$20$&$^{\rm o}{\rm C}$\\
		$H$                &mean ocean depth&$4500$& ${\rm m}$\\
		$V$ &control volume& $8250\times 4.5\times 300$& ${\rm Km}^{3}$\\
		\hline
	\end{tabular}
\end{table}

\begin{figure}[h]
	\centering
	\includegraphics[width=.6\textwidth]{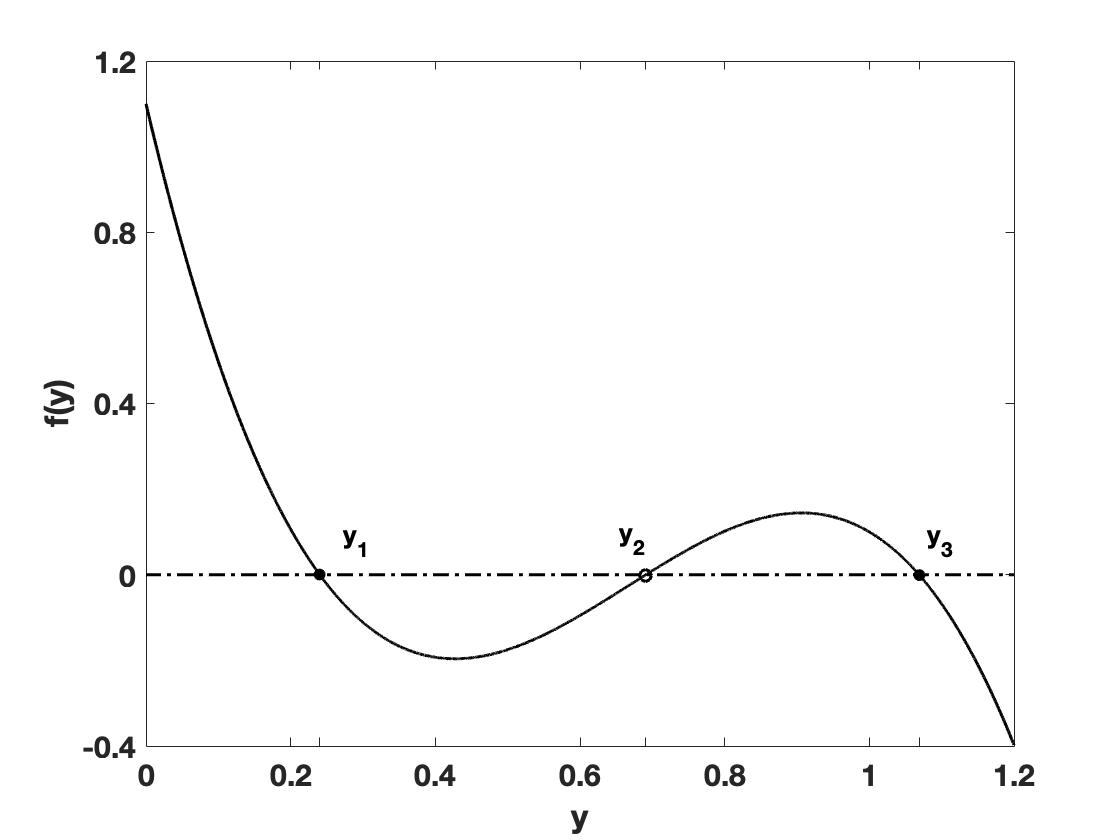}
	\caption{Equilibrium states for deterministic system about salinity difference (\ref{e11}) with freshwater flux $\bar{F}=1.1$ and time scale ratio $\mu^2=6.2$.}
	\label{fig1}
\end{figure}
%

The system expressed by (\ref{e8}) is a slow-fast dynamical system with $x$ being the fast variable  and $y$   the slow variable. Here, only the dynamics of the slow variable $y$ is of interest. It is reasonable to assume that the temperature difference $x$ in the high latitudes is constant for a large $\alpha$ \cite{Dijkstra2013}. Therefore, the temperature is clamped to the prescribed value $x=1+\mathcal{O} (\alpha ^{-1})$. Consequently, the dynamics for slow variable $y$, the salinity difference between equatorial and pole regions,  can be expressed as
\begin{align}\label{e11}
	\dfrac{{\rm d} y}{{\rm d}t } = f(y),
\end{align}
where $f(y)=\bar{F}-y\left(1+\mu^2 \left(1-y\right)^2\right)$.

The equilibrium states are obtained by  setting the right hand side of equation \eqref{e11} equal to zero.  With freshwater flux $\bar{F}=1.1$ and time scale ratio $\mu^2=6.2$, values for the three equilibrium states are  $y_1=0.2402$, $y_2=0.6911$ and $y_3=1.0687$, respectively, as shown in Finger 1.


Now   suppose that the freshwater flux $\bar{F}$ has a stochastic component, and consider a stochastic thermohaline circulation  model under additive Gaussian noise,
\begin{align} \label{e12}
\begin{split}
{{\rm d} Y(t)} =\left(\bar{F}-Y(t)\left(1+\mu^2 \left(1-Y(t)\right)^2 \right)\right){\rm d}t+\epsilon {\rm d}B(t) ,~~~~0< t < T,
\end{split}
\end{align}
where $\epsilon $ is noise intensity, $T$ is transition time and $B(t) $ is a standard Brownian motion. This model \eqref{e12} has been widely used to study the thermohaline circulation under stochastic fluctuations of  salinity forcing \cite{Velez2001}. Without stochastic fluctuations, i.e. $\epsilon=0$, the salinity difference $y$ always reaches one of the two stable equilibria located at the solid points as shown in Figure 1.

Since this paper is mainly interested in transitions between the two equilibrium states $y_1$ and $y_3$,  (\ref{e12}) is investigated under the following bridge condition,
\begin{align}\label{eq13}
Y(0)=y_1, ~~~~Y(T)=y_{3}.
\end{align}

\section{Transition of the maximum likelihood state}

In this section, we shall study how the maximum likelihood state of (\ref{e12}) under the bridge condition \eqref{eq13}  evolves with time. Here the maximum likelihood state at time $t$ for (\ref{e12}) is defined as
\begin{align} \label{eq14}
	\psi(t)=\underset{y \in \mathbb{R}}{{\rm arg~ max}} ~p(y,t|y_3,T; y_1,0),
\end{align}
where  $p(y,t|y_3,T; y_1,0)$ is the density of $Y(t)$ in (\ref{e12})  at $Y(t)=y$ given $Y(0)=y_1$ and $Y(T)=y_3$. Note that $\psi(t)$ is the most probable state at time $t$ among all the possible states given the initial and final states.

Let $p(y,t|y_1,0)$ be the density of $Y(t)$ at  $Y(t)=y$ conditioned on $Y(0)=y_1$.  It satisfies the following forward Fokker-Planck equation \cite{Duan2015,Risken1996},
\begin{align}\label{e15}
\begin{cases}
\dfrac{\partial}{\partial t} p(y, t|y_1,0)=-\dfrac{\partial}{\partial y} \left(\left(\bar{F}-y\left(1+\mu^2 \left(1-y\right)^2\right) \right)p(y,t|y_1,t_1)\right)+\dfrac{\epsilon ^2}{2} p(y,t|y_1,t_1), ~~~{\rm for}~~~ t>0,\\
~~~p(y,0|y_1,0)=\delta(y-y_1).
\end{cases}
\end{align}

Let $p(y_3,T|y,t)$ be the density of $Y(t)$ at  $Y(t)=y$ conditioned on $Y(T)=y_3$.  It satisfies the following backward Fokker-Planck equation \cite{Risken1996},
\begin{align}\label{e16}
\begin{cases}
\dfrac{\partial}{\partial t} p(y_3, T|y,t)=-\left( \bar{F}-y\left(1+\mu^2 \left(1-y\right)^2\right)\right)\dfrac{\partial}{\partial y} p(y,t|y_1,t_1)-\dfrac{\epsilon ^2}{2} \dfrac{\partial ^2}{\partial y^2}p(y_2,T|y,t),~~~{\rm for~~~ t<T},\\
~~~p(y_3,T|y,T)=\delta(y_3-y).
\end{cases}
\end{align}

Note that
\begin{align} \label{eq17}
	p(y, t |y_3,T; y_1, 0) =\dfrac{p(y_2, T|y,t) p(y,t|y_1,0)}{p(y_3,T|y_1,0)},
\end{align}
 and  $p(y_3,T|y_1,0)$ is a constant independent of $y$ and $t$, it follows from \eqref{eq14} and \eqref{eq17} that
\begin{align}\label{e18}
\psi(t)=\underset{y \in \mathbb{R}}{{\rm arg~ max}} ~p(y_3, T|y,t) p(y,t|y_1,0).
\end{align}

\begin{figure}[h]
	\centering
	\subfigure[]{
		\includegraphics[width=.48\textwidth]{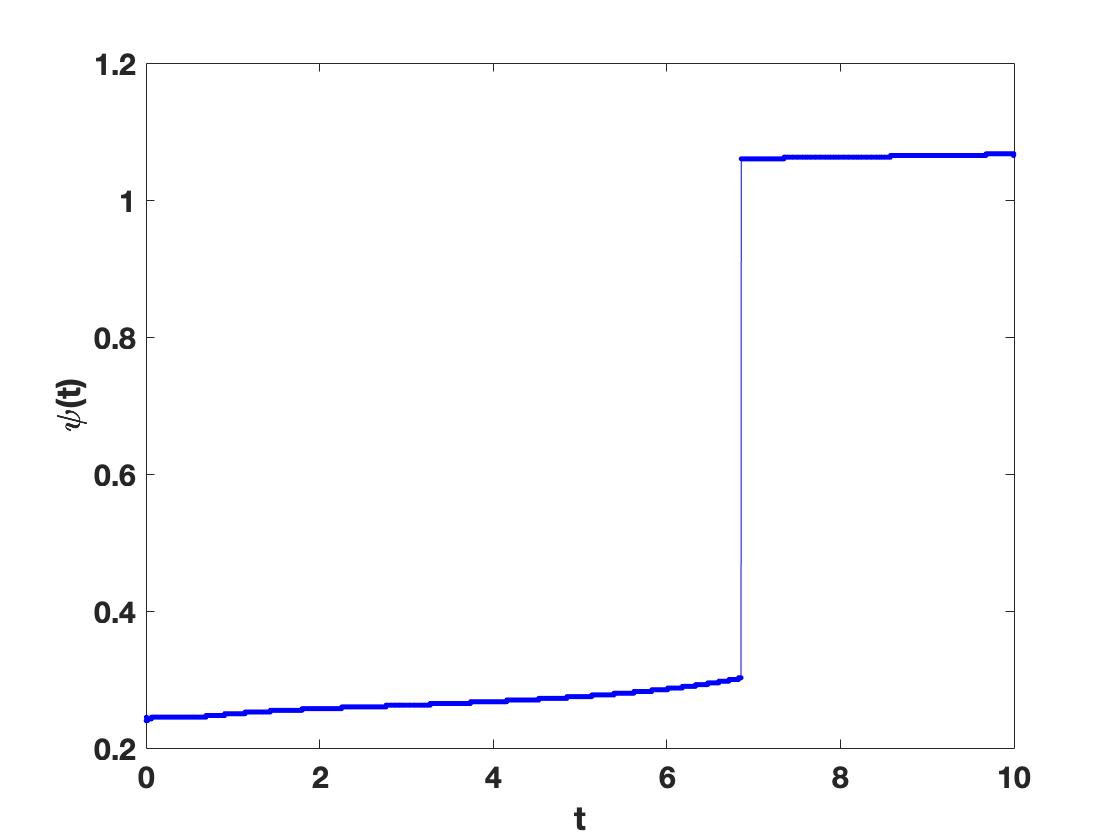}
	}
	\subfigure[]{
		\includegraphics[width=.48\textwidth]{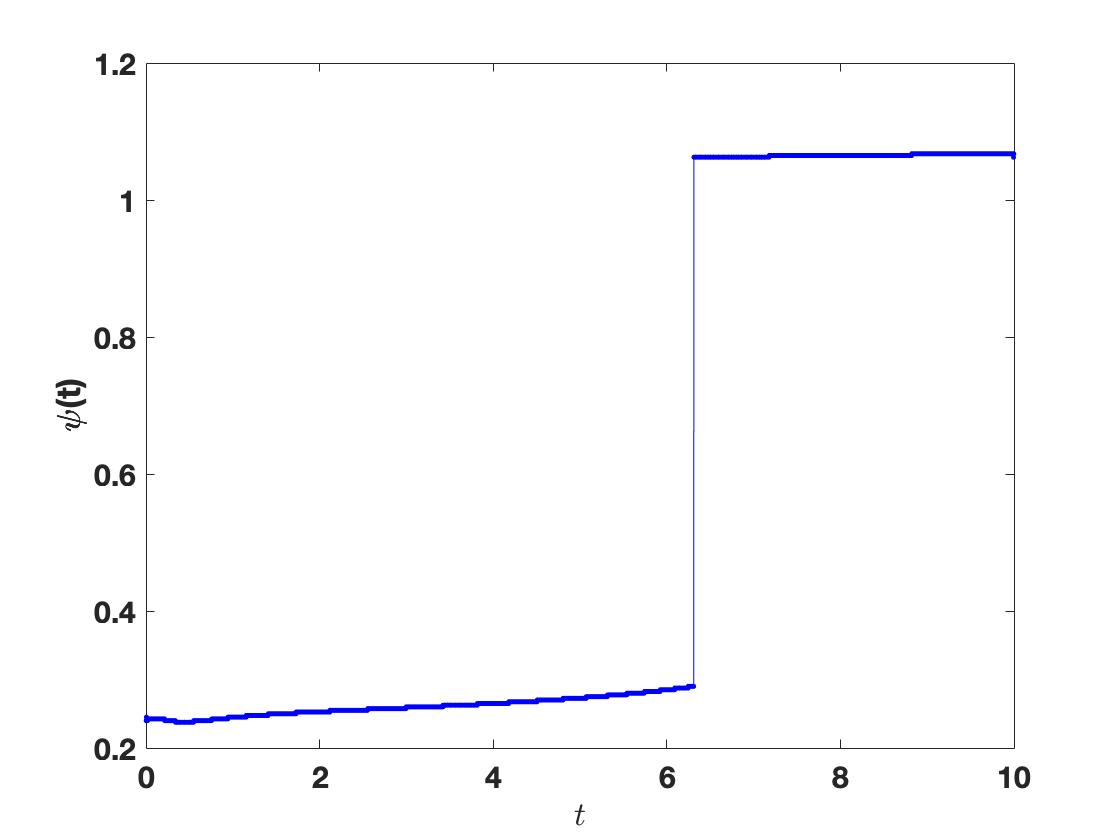}
	}
	\caption{The maximal likelihood transition pathway of stochastic system for salinity difference (\ref{e12}) under condition \eqref{eq13} with freshwater flux $\bar{F}=1.1$, time scale ration $\mu^{2}=6.2$ and noise intensity (a) $\epsilon = 0.20$, (b) $\epsilon =0.25$. }
	\label{fig3}
\end{figure}

The values of $p(y,t|y_1,0)$ and $p(y_3,T|y,t)$  in (\ref{e18}) for every $y$ and $t$ can be obtained by numerically solving (\ref{e15}) and (\ref{e16}), respectively.  For detailed discuss on the numerical methods for (\ref{e15}) and (\ref{e16}), see \cite{YangSun2020}.

For parameters $\bar{F}=1.1$, $\mu^{2}=6.2$ and noise intensity $\epsilon=0.20, 0.25$, the graphs of $\phi(t)$ are shown in Figure \ref{fig3} (a) and (b), respectively. It can be seen from the Figure 2 (a) that a jump occurs at $t\approx 6.86$, during the state transition, and the state changes rather slowly before and after the jump, which indicates an abrupt change. As shown in Figure 2 (b), the jump occurs at $t\approx 6.32$.  Our simulation shows that the time for the abrupt change is greatly influenced by the strength of the noise. Figure 3 shows how the time of abrupt change varies with noise intensity $\epsilon$.  As can be seen from the figure, the abrupt change occurs earlier with the stronger noise.

\begin{figure}[h]
	\centering
	\includegraphics[width=.6\textwidth]{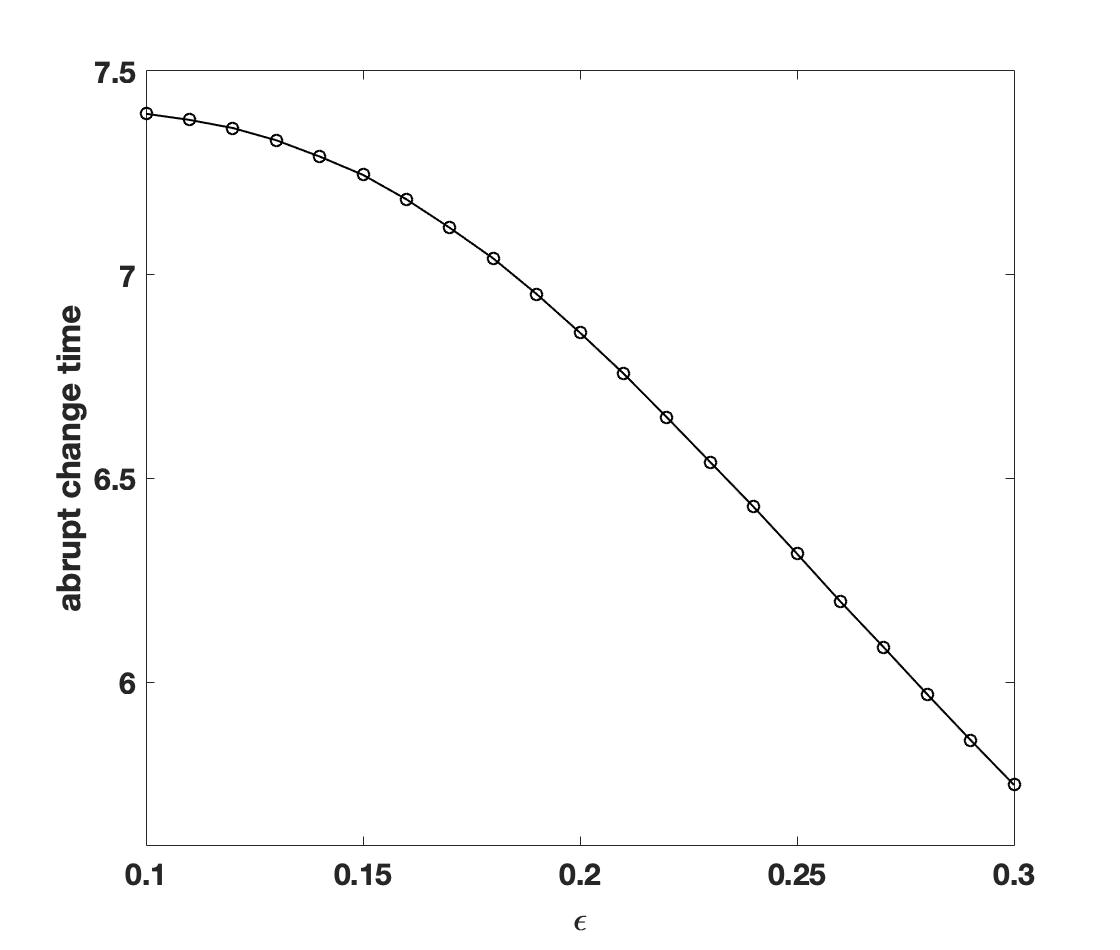}
	\caption{The time for abrupt change $t$ as function of noise intensity $\epsilon$ for stochastic system of salinity difference  \eqref{e12} under the condition $Y(0) = 0.2402,$ and $Y(10) = 1.0687$.}
	\label{fig6}
\end{figure}

\section{Conclusion}

Based on a simplified stochastic thermohaline circulation system, we  study how the maximum likelihood state evolves during the transition from one equilibrium state to the other. It is shown, with the help of simulation results, that there is a jump in the transition pathway of the maximum likelihood state. The jump, which indicates an abrupt change for the salinity variability in the thermohaline circulation system, is helpful to understand the abrupt climate change as revealed by various climate records.

\section*{Acknowledgements}
This work is supported by the National Natural Science Foundation of China (grant No.11801192, 11531006 and 11771449).

\section*{Data Availability Statement}
The data that support the findings of this study are openly available on GitHub \cite{code}.

\bibliographystyle{unsrt}
\bibliography{maxlikehood}

\begin{thebibliography}{10}

\bibitem{Weaver1992}
A.~J. Weaver, J.~Marotzke, P.~F. Cummins, and E.~S. Sarachik.
\newblock {S}tability and {V}ariability of the {T}hermohaline {C}irculation.
\newblock {\em Journal of Physical Oceanography}, 23(1):39--60, 1993.

\bibitem{Monahan2002}
A.~H. Monahan.
\newblock {S}tabilization of {C}limate {R}egimes by {N}oise in a {S}imple
  {M}odel of the {T}hermohaline {C}irculation.
\newblock {\em Journal of Physical Oceanography}, 32(7):2072--2085, 2002.

\bibitem{Velez2001}
P.~V{\'e}lez-Belch{\'i}, A.~Alvarez, P.~Colet, J.~Tintor{\'e}, and R.~L. Haney.
\newblock {S}tochastic {R}esonance in the {T}hermohaline {C}irculation.
\newblock {\em Geophysical Research Letters}, 28(10):2053--2056, 2001.

\bibitem{Cessi1994}
P.~Cessi.
\newblock A {S}ilmple {B}ox {M}odel of {S}tochastically {F}orced {T}hermohaline
  {F}low.
\newblock {\em Journal of Physical Oceanography}, 24(9):1911--1920, 1994.

\bibitem{Timmermann2000}
A.~Timmermann and G.~Lohmann.
\newblock {N}oise-induced {T}ransitions in a {S}implified {M}odel of the
  {T}hermohaline {C}irculation.
\newblock {\em Journal of Physical Oceanography}, 30(8):1891--1900, 2000.

\bibitem{Lorenzo2008}
M.~N. Lorenzo, J.~J. Taboada, I.~Iglesias, and I.~Alvarez.
\newblock {T}he {R}ole of {S}tochastic {F}orcing on the {B}ehavior of
  {T}hermohaline {C}irculation.
\newblock {\em Trends and Directions in Climate Resesrch}, 1146:60--86, 2008.

\bibitem{Knight2005}
J.~R. Knight, R.~J. Allan, C.~K. Folland, M.~Vellinga, and M.~E. Mann.
\newblock {A} {S}ignature of {P}ersistent {N}atural {T}hermohaline
  {C}irculation {C}ycles in {O}bserved {C}limate.
\newblock {\em Geophysical Research Letters}, 32(20):L20708, 2005.

\bibitem{Zhang2010}
S.~Zhang, A.~Rosati, and T.~Delworth.
\newblock {T}he {A}dequacy of {O}bserving {S}ystems in {M}onitoring the
  {A}tlantic {M}eridional {O}verturning {C}irculation and {N}orth {A}tlantic
  {C}limate.
\newblock {\em Journal of Climate}, 23(19):5311--5324, 2010.

\bibitem{Averyanova2017}
E.~A. Averyanova, A.~B. Polonsky, and V.~F. Sannikov.
\newblock {T}hermohaline {C}irculation in the {N}orth {A}tlantic and its
  {S}imulation with a {B}ox {M}odel.
\newblock {\em Izvestiya Atmospheric and Oceanic Physics}, 53(3):359--366,
  2017.

\bibitem{Rahmstorf2002}
S.~Rahmstorf.
\newblock {O}cean {C}irculation and {C}limate {D}uring the {P}ast 120,000
  {Y}ears.
\newblock {\em Nature}, 419(6903):207--214, 2002.

\bibitem{Dijkstra2013}
H.~A. Dijkstra.
\newblock {\em Nonlinear Climate Dynamics}.
\newblock Cambridge University Press, 2013.

\bibitem{Duan2015}
J.~Duan.
\newblock {\em An Introduction to Stochastic Dynamics}.
\newblock Cambridge University Press, 2015.

\bibitem{Risken1996}
H.~Risken and T.~Frank.
\newblock {\em The Fokker-Planck equations: Methods of Solution and
  Applications}.
\newblock Springer, 1996.

\bibitem{YangSun2020}
F.~Yang and X.~Sun.
\newblock {O}n {D}ynamics of the {M}aximum {L}ikelihood {S}tates in
  {N}onequilibrium {S}ystems.
\newblock {\em Journal of Statistical Physics}, 181(3):753--760, 2020.

\bibitem{code}
F.~Yang.
\newblock Code.
\newblock Github, 2020.
\newblock
  \url{https://github.com/yangfang0914/On-the-abrupt-change-of-the-maximum-likelihood-state-in-a-simplified-stochastic-thermohaline-circula.git}.

\end{thebibliography}

%
\end{document}